\documentstyle{amsppt}
\headline={\hfil}
\footline={\ifnum\pageno<2 \hfil \else \hss\folio\hss \fi}
\output={\plainoutput}
\hcorrection{.65truein}
\comment
\magnification=\magstep1

\hsize=6.5truein
\vsize=9truein
\abovedisplayskip=15pt plus 3pt minus 3pt
\abovedisplayshortskip=10pt plus 3pt
\belowdisplayskip=15pt plus 3pt minus 3pt
\belowdisplayshortskip=10pt plus 3pt minus 2pt
\endcomment

\def\qed{\hfill $\square$}
\def\os{\overset}
\def\us{\underset}

\parindent = 20pt      

\textfont\ttfam=\tentt \def\tt{\fam\ttfam\tentt}%

\def\q{\quad} 

\def\vs{\vskip.3cm}
\def\vsk{\vskip.5cm}
\def\vvs{\vspace{2\jot}}
\def\noi{\noindent}
\def\npb{\nopagebreak}

\def\idealP{{\text{\tt P}}}
\def\idealq{{\text{\tt q}}}
\def\q{{\text{\tt q}}}

\def\ideals{{\text{\tt s}}}
\def\idealc{{\text{\tt c}}}

\def\Spec{\text{\rm Spec}}

\def\Het{{\text{\rm H}}_{\text{\rm \'et}}}
\def\Hzar{{\text{\rm H}}_{\text{\rm Zar}}}

\def\exp{\text{\rm exp}}
\def\ad{\text{\rm ad}}
\def\Pic{\text{\rm Pic}}
\def\kalg{k\text{-alg}}
\def\IN{\Bbb N}
\def\Gm{\text{\pmb G}_{m}}

\def\J{\Cal J}

\def\fG{\frak G}
\def\fL{\frak L}
\def\fU{\frak U}
\def\fX{\frak X}
\def\fY{\frak Y}
\def\X{\frak X}
\def\fXred{\frak X_{\text{\rm red}}}
\def\Xred{\frak X_{\text{\rm red}}}
\def\fT{\frak T}
\def\fg{\frak g}
\def\fgl{\frak{gl}}

\def\p{{\text{\tt p}}}
\def\P{{\text{\tt P}}}
\def\Q{{\text{\tt Q}}}

\def\fh{\frak h}
\def\fk{\frak k}
\def\fr{\frak r}
\def\fa{\frak a}
\def\kdiag{k\text{-diagonalizable}}
\centerline{}
\vs
\centerline{\bf   LOCALLY TRIVIAL PRINCIPAL HOMOGENEOUS SPACES AND}
\vs
\centerline{\bf   CONJUGACY THEOREMS FOR  LIE ALGEBRAS}

\vskip.5in
\centerline{\smc  A. Pianzola
\footnote{Supported by NSERC operating grant A9343.
\newline
\indent 1991 {\it Mathematics Subject Classification:} 14L15, 17B01 and
22E65.}  }

\vs
\centerline{Department of Mathematical Sciences}
\centerline{University of Alberta}
\centerline{Edmonton, Alberta}
\centerline{Canada  T6G 2G1}

\topmatter
\abstract\nofrills{}
We link locally trivial principal homogeneous spaces over $\text{Spec}\,R$ to the question of
conjugacy of maximal abelian diagonalizable subalgebras of $\fg \otimes 
R.$
\endabstract
\endtopmatter

\baselineskip = 18pt 

  Throughout  $k$ will denote a field of characteristic zero. Unless specifically mentioned
otherwise {\it all  algebras, tensor products, vector spaces, and schemes  are over }
$k.$ 
\vs

Let $\fg$ be a finite dimensional split semisimple
Lie algebra over  $k$ . Of central importance to
classical Lie theory is Chevalley's theorem asserting the conjugacy of all split
Cartan subalgebras of $\fg$. We are interested in an analogue of this result for
Lie algebras of the form $\fg(R) := \fg \otimes R,$  where $R$ is an associative commutative
unital
$k\text{-algebras}$ (Recall that $\fg(R)$ is viewed {\it as an algebra
over} $k. $ In general, these algebras are infinite dimensional.) A well understood  example 
is the  case of  the algebra
$R = k[t,t^{-1}]$ of Laurent polynomials. Then
$\fg(R)$ is the so called loop algebra of $\fg$ that one encounters on the
realizations of  non-twisted affine Kac-Moody Lie
algebras.  In this case the appropriate version of
conjugacy  is  due to Peterson and Kac (see Remark 2(iii) below).

Let $\fh$ be a split Cartan subalgebra of $\fg.$  Then $\fh \simeq \fh
\otimes 1$ is {\it not} in general a Cartan subalgebra of $\fg(R)$ (since it is
not self normalized unless $R = k$).  The split Cartan subalgebras of $\fg$ 
are examples of {\it abelian}  $k${\it -diagonalizable subalgebras} 
of $\fg(R)$, 
namely   of subalgebras $\fa$ of $\fg(R)$ such  that 
\vs

\item {\rm(i)} $\fa $ is abelian.

\item{\rm(ii)} All elements of $\fa$ are $ k\text{-diagonalizable}:$ 
If $\p$ belongs to $\fa$ then
$\text{ad}_{\fg (R)}\p ,$ when viewed as a $k$-linear endomorphism of $ \fg(R),$
is diagonalizable.

\vs
\noi (Any subalgebra of 
$\fg(R)$ satisfying (ii) is abelian, but no harm is
done by emphasizing this last). If in addition 

\vs
\item {\rm(iii)} No subalgebra of $\fg (R)$ satisfying (i) and (ii) above
properly contains $\fa.$
\vs
 \noi then $\fa$ is a {\it maximal abelian} $k\text{\it -diagonalizable}$ {\it subalgebras}, or MADs for 
short. (We will see later that  split Cartans of $\fg$ are MADs of 
$\fg(R)$ if and only if $\Spec(R)$ is connected).

 Since these type of subalgebras play a crucial role in
understanding
$\fg(R)$ and its representations in both the  finite dimensional and affine Kac-Moody case, it
is natural and relevant to ask if all MADs of
$\fg(R)$ are conjugate under some suitable subgroup of  $\text{Aut}_{k\text{-Lie}}\,\fg 
(R).$ The natural choice for this subgroup (because of functoriality on 
$R$ and compatibility with the usual results in the case
of a base field), is the 
group $\fG(R)$ of $R$-points of the 
corresponding simply connected 
Chevalley-Demazure group, acting on $\fg(R)$ via the adjoint 
representation. As we shall see, the answer to this question is quite 
interesting and  related to the triviality of certain principal 
homogeneous spaces over $\Spec(R)$. 

Again  by analogy with the finite dimensional case, one expects regular elements to 
play a special role in the 
problem at hand. The correct functorial definition for these elements 
is as follows \footnotemark"$^{2}$". Let $f_{\text{reg}} \in  S(\fg^*)$ be the polynomial
function defining the basic Zariski open dense set of regular elements 
of $\fg$ (see [Bbk2] Ch.
VII). Since 
$f_{\text{reg}}$ is defined over
$k,$ we can think of it as a polynomial function on the free $R\text{-module}$
$\fg(R).$ An element $\p$  of
$\fg(R)$ will be said to be {\it regular} \;if $f_{\text{reg}}(\p)$ is a
unit of $R.$ Finally, a MAD is said to be {\it regular} if it 
contains a regular element.
$\footnotetext"$^{2}$ "{ J-P Serre suggested this definition. See also
Expose XIV of SGA3.}$

Here then is our main result.
\vs

\proclaim{Theorem 1}Let $\fg$ be  a finite dimensional
split semisimple Lie algebra over $k,$ and  $\fG $ its simply connected
Chevalley-Demazure group scheme.  Let $\fX = \Spec(R)$ be a  
connected affine scheme and $\fXred$ the corresponding reduced scheme. Then
\vs
\item {\rm (i)} If $\fa$ is an abelian $k\text{-diagonalizable}$
subalgebra of $\fg(R)$ then $\text{dim}_k(\fa) \leq \text{rank}(\fg).$ If this
is an equality then $\fa$ is maximal.

\item {\rm (ii)} Assume that $\fX(k) \neq \emptyset$. 

{\rm (a) (Regular conjugacy).} If the Picard group of $\fXred$ is trivial 
then all regular maximal
abelian
$k\text{-diagonalizable}$  subalgebras of  
$\fg(R)$ are  conjugate under
$\fG(R)$.

{\rm (b) (Full conjugacy).} Consider the following property on $\fX$. 

\noi {\bf (TLT)} \text{\rm 
(Triviality of locally trivial Levi torsors):} If $\fL$ is the Levi subgroup 
of a standard parabolic subgroup of $\fG,$ then any locally 
trivial principal homogeneous
space for $\fL$ over $\fXred$ is
trivial.

\noi If {\rm (TLT)} holds, then all maximal abelian
$k\text{-diagonalizable}$ subalgebras of
$\fg(R)$ are regular (and hence all conjugate by \text{ \rm(a)}).

\endproclaim

 The main idea behind  the proof of Theorem 1 is to evaluate the different primes
of $\fX$ at a given  $\kdiag$ element of $\fg(R).$ Each of these
evaluations puts us in the finite dimensional case where conjugacy is known to
hold. One then is forced to look at assumptions on $\fX$  that allow
all of these finite dimensional conjugacies to be ``pasted together" to create an element of
$\fG(R).$ The proof of the first part of Theorem 1 is straightforward 
and is given earlier in the paper after developing some basic properties of 
$\kdiag$ elements. This is followed by a series of results that 
conclude in Proposition 11 with the translation of the conjugacy question to one on the 
triviality of certain torsors (= principal homogeneous spaces. See 
Remark 2(ii) below) over $\Spec (R)$. An 
induction argument is then used to prove the second part of 
the main Theorem. The paper concludes with an interesting 
example. 
\vs

\noi{\bf 2 Remarks.}  

 (i) Most of the assumptions  of the Theorem
are not  superflous. 
There exist  rings $R$ leading to non regular
MADs, and if $\Pic(\fX) \neq 0$ regular MADs need not be conjugate. 
The connected assumption on $\fX$ is needed in part (i) of the 
Theorem but is not crucial 
for part (ii) (which holds if $\fX$ has a finite number of connected components 
each of which satisfies the assumptions of the Theorem). On the other hand the assumption on the 
existence of a rational point, namely of a 
maximal ideal $x_{0}$ such that $R/x_{0} \simeq k$, is central to 
the proof. 

(ii) (See [Mln] and [DG] for details). Let $\fX$ be a $k$-scheme, and 
let $\fL$ be an algebraic $k$-group. A (right) $\fL$-torsor over $\fX$ (also 
called an $\fX$-torsor under $\fL$) is an $\fX$-scheme $\fY$ on which 
 $\fL_{\fX} := \fL \times_{\Spec (k)} \fX$ acts on the right, and which is
locally isomorphic  to $\fL_{\fX}$ for  the flat topology of $\fX$ (with 
$\fL_{\fX}$ 
acting on itself by right
multiplication). Thus there 
exist flat and locally finitely presented morphisms $\phi_{i}: \fU_{i} 
\rightarrow \fX$  with $\fX = \cup \phi_{i}(\fU_{i})$ and $\fY \times_{\fX} 
\fU_{i} \simeq \fL _{\fU_{i}} :=\fL_{\fX} \times_{\fX} 
\fU_{i}$ (these isomorphisms preserving the respective $\fL_{\fU_{i}}$ 
actions). If our group is smooth the $\phi_{i}$ may be taken to 
be \'etale, and then just as with principal bundles in differential geometry 
(of which torsors are a suitable algebraic analogues) we can attach to 
the isomorphism class of a torsor $\fY$ as above an element of 
$\Het^{1}(\fX,\fL_{\fX})$ (C\v ech cohomolgy on the \'etale site of $\fX$ 
with coefficients on the group sheaf $\fL_{\fX}$). This is an 
injective procedure, and it is 
surjective if $\fX$ is affine.  $\Het^{1}(\fX,\fL_{\fX})$ is a set with a distinguished element, 
namely the isomorphisms class of the trivial torsor $\fL_{\fX}$ acting 
on itself by right multiplication. 
If the $\phi_{i}$ can be taken to be  open immersions (the $\fU_{i}$ 
can then be thought as an honest open cover of $\fX$ in the Zariski topology), the 
torsor $\fY$ is said to be {\it locally trivial}. Their isomorphism 
classes are then parametrized (again assuming $\fX$ affine and $\fL$ 
smooth) by 
$\Hzar^{1}(\fX,\fL_{\fX})$

(iii) Condition TLT   varies with the type of $\fg$ and the
nature of the base field $k.$ For  the condition to hold for all 
types it suffices to assume that $\fX$ has the following property.
\medskip
\noi {\bf (TRT)} \text{\rm 
(Triviality of locally trivial reductive torsors):} { \it If $\fL$ is 
a (connected) 
split reductive  $k$-group then any locally 
trivial principal homogeneous
space for $\fL$ over $\fXred$ is
trivial.}
\medskip
There are two important examples of rings with this property, namely those $R$ which modulo their 
nilradical equal either

	\item{(a)} $ k[t_1,...,t_n]$, or
	\item{(b)} $ k[t_{1}^{\pm 1},...,t_{n}^{\pm 1}]$

Case (a) follows fron the work of Raghunathan and Ramanathan, and of
Raghunathan  dealing with the triviality of certain 
torsors over  algebraic affine space
([Rgn1] Theorem 2.2. See also [CTO]). Case (b) reduces to case (a). For $n = 1$ this is easy 
since every locally trivial principal 
homogeneous space under $\fG$
over the punctored affine line, 
 extends to one over the  whole affine line (In fact because our 
 torsors are locally trivial, one can directly show that (TLT) holds 
 for $\fX = \Spec(k[t^{\pm 1}])$ by means of a standard argument. See for 
 example the proof of proposition 3(i) of [Pzl2]).
Note that this recovers the conjugacy theorem of Peterson-Kac in the case of
untwisted affine Kac-Moody Lie algebras (see [PK]). 
The general case follows by an induction argument due to Gille [Gll]. 
Ccase (b)  is of great importance because of 
toroidal Lie algebras (see [Pzl3]). 

The proof of the  Theorem is slightly easier if one assumes (TRT) 
rather that (TLT) (see the Remark following Proposition 9). Non 
standard examples where (TLT) holds
can be
found with the aid of Th\'eor\`eme 6.13 of [CTS]. 

(iv) Note that we are dealing with the triviality of certain {\it 
algebraic} principal bundles over the {\it global} space $\fX = \Spec (R).$
In particular one is not allowed to replace $R$ by any of its 
localizations or completions (indeed, conjugacy may hold for all 
localization of $R$ yet fail for $R$ itself). That we are in the 
algebraic setup forces us to work, even if the base field $k$ is $\Bbb R$ 
or $\Bbb C$, with the Zariski topology and the complications that this 
entails for fibrations (compare for example the triviality of vector bundles over 
affine space in the classical case by contractability, with its 
``Serre's conjecture'' 
algebraic counterpart: Theorems of Quillen and Suslin). The work of 
Raghunathan is here crucial.

It is important to observe that though MADs behave somehow 
functorialy on $R$ (Lemma 5 and Proposition 6), MADs are not 
$R$-modules. The point is that $k$-diagonalizability is lost in 
general by scaling under elements of $R$ which are not in $k$. In fact 
the role of $k$-diagonalizability is crucial but deceivingly subtle 
and may at times be easily overlooked.

\vs

\noi {\bf 3 Notation and conventions.}

Throughout   $\fg$ and  $\fG$ will be as in the statement 
of Theorem 1. The category of 
commutative associative unital 
$k\text{-algebras}$ will be denoted by $k\text{-alg}$. If $R$ is in 
$\kalg$  
the residue field of 
an element $x$ of $\Spec(R) =\fX$ will be denoted by $k(x).$ For convenience in what follows the group 
$\fG(k(x))$ will be denoted simply by
$\fG(x),$ and the corresponding  
 group  homomorphism $ \fG(R)\to \fG(R/x) \subset \fG(x)$  by $\idealP \mapsto \idealP(x).$

The constructions of the last paragraph  can be repeated, {\it mutatis
mutandi}\, , if we replace
$\fG$ by its Lie algebra functor
$\fg(-).$ Since $\fg$ is finite dimensional we have 
$\fg(-) = \text{Hom}_{k\text{--alg}}
(S(\fg^*) , -).$ Thus $\fg(S) = \fg \otimes S$ for any $S$ in $k\text{-alg}$. In
particular $\fg \simeq \fg(k).$

Along similar lines if $V$ is a vector space over $k$, $S$ is in $k\text{-alg}$, 
and $x \in \fX$; we will denote $V \otimes S$ by $V(S)$ and  $V(k(x))$
by $V(x).$

Let $k[\fG]$ be the coordinate ring of $\fG.$ There is then a dual nature to $\fG.$ It can be 
thought as the 
scheme $\text{Spec}(k[\fG])$ or as the functor $\text{Hom}_k(k[\fG], \,- \,)$ 
from $k\text{-alg}$ into the category of groups. We shall make use of both  these
manifestations. The following example may help clarify these  ideas.
\vs
\noi{\bf Example. } Take $\fg$ to be of type $\frak{sl}_n$, $k = 
\Bbb{R}$, and $R =
\Bbb{R}[t,t^{-1}].$ Then $\fG = \text{\bf  SL}_n = \text{Hom}_{\Bbb{R}\text{-alg}}(\Bbb{R}[\fG], \,-
\,)$ with
$\Bbb{R}[\fG] = \Bbb{R}[x_{ij}]/(det - 1).$
$\fG(R)$ (respectively $\fg(R)$) is the group (Lie algebra) of $n \times n$
matrices of determinant 1 (trace 0) with entries in the ring $R$.  We have
$\fX = \{<f(t)>, <g(t)>, \{0\} \}$ where $f(t)$ and $g(t)$ are irreducible of
degree 1 and 2 respectively, and $f(0) \neq 0.$
The corresponding residue fields of these three types
of primes are isomorphic to
$\Bbb{R}$, $\Bbb{C},$ and $\Bbb{R}(t)$ respectively.  If $\P$ is an element 
of $\fG(R)$ and
$x \in \fX,$ then $\P(x)$ is simply the matrix obtained by reducing
$\text{mod}\,x$ each of the entries of $\P.$ Similarly for $\p \in 
\fg(R).$ 
\vs

\vs
 We begin by  recalling an important fact of which we will make repeated use
in what follows.
\proclaim{Proposition 4}
Let $\fa$ be an abelian $\kdiag$ subalgebra of $\fg(R).$ Assume $R$ is 
an integral domain, and let $K$ denote its field of quotients. Then there
exists a split Cartan subalgebra $\frak k$ of the $K$-Lie algebra $\fg (K)$ with $\fa \subset
\frak k.$
\endproclaim 
\demo {Proof}
 See Seligman [Slg]. See also [Bbk2] Ch.~8 Exercise \S3.10(b)  \qed
\enddemo

\vs
The following result is straightforward.

\vs
\proclaim {Lemma 5}  Let $x\in \fX$. 

\item {\rm (i)} If $\fa$ is an abelian $k\text{-diagonalizable}$ subalgebra of 
$\fg(R)$ then  
$\fa(x) := \{ \p(x) : \p \in \fa\}$ is 
abelian and $k\text{-diagonalizable}$ when viewed as a  subalgebra of either $\fg(R/x)$ or $\fg(x).$  In 
particular, if $\p$ is $\kdiag$ then so 
is $\p(x).$
\item {\rm (ii)} If $\p \in \fg(R)$ is regular then  $\p (x)$ is regular.\qed
\endproclaim

\vs

Let $\rho:\fg \to \fgl (V)$ be a finite dimensional representation of $\fg.$
 For any $S$ in $k\text{-alg}$ we then get a representation
$
\rho_S :\fg (S)\to \fgl \big(V(S)\big)
$
of the $S\text{-Lie}$ algebra $\fg (S).$ If $\p \in \fg(S)$ then 
$\rho_S(\p)^ m$ is an $S$-linear endomorphism of 
the  free 
 $S$-module of finite rank $V(S).$ It is meaningful therefore to 
 consider its trace.
\vs
\proclaim {Proposition 6} Let $\rho:\fg \to \fgl (V)$ be a finite dimensional
representation of 
$\fg.$ Assume $\fX$ is connected and reduced. If $\p\in \fg (R)$
is $k\text{-diagonalizable}$ then. 

\item {\rm (i)} $\text{\rm Tr}\, \rho_R(\p)^ m\in k$ for all $m\in \IN.$
\item {\rm (ii)} If $\p(y) = 0$ for some $y \in \fX$ then $\p = 0.$
\endproclaim 

\demo {Proof} We reason in stages.

\noi Step 1: {\it Reduction to the noetherian case.}  Let $R'$ be a finitely generated subalgebra of $R$ such 
that $\p$ can be viewed as an element $\p'$ of $\fg(R') \subset 
\fg(R).$ Clearly $\fX' = \Spec(R')$ is connected and reduced, and 
$\p'$ is a $\kdiag$ element of $\fg(R')$. Since $\text{\rm Tr}\, \rho_{R'}(\p')^ m = \text{\rm Tr}\, 
\rho_R(\p)^ m$ and $\p(y) = 0$ implies $\p'(y \cap R') = 0$, it will 
suffice to establish the result under the assumption that $R$ is 
noetherian.

\noi Step 2:  $\fX$ {\it integral.}  Let $K$ be the field of quotients of $R$. Since
$\text{ad}_{\fg(K)}\p$ is semisimple $\rho_K(\p)$ acts semisimply on 
$V(K).$  We claim that the eigenvalues of $\rho_K(\p)$ (in the algebraic
closure of $K)$  belong to $k$. To see this
put $\p$ inside a split Cartan $\frak k$ of $\fg(K)$  as in Proposition 
4,
and fix a base $\Pi = \{\alpha _1\dots \alpha _\ell\}$  of the corresponding root
system $\Delta =\Delta (\fg(K),\frak k).$ If
$\omega ^\vee_1,\dots,\omega ^ \vee_\ell$ are the fundamental coweights corresponding to
the $\alpha _i\text{'s}$ then $\frak k =\oplus K\omega ^ \vee_i.$
Now if  we write $\p=\sum c_i\omega ^\vee_i$
then the $c_i\text{'s}$ are eigenvalues of
$\text{ad}_{\fg(K)}\p$ and therefore belong to $k$. Since the weight space decomposition with 
respect to $\fk$ of the 
the $\fg(K)$-module afforded by $\rho_K(\p)$ are rational linear 
combinations of the  $\alpha _i\text{'s}$ the claim follows. 
From the above we conclude that  $\text{\rm Tr}\, \rho_K(\p)^ m :=\lambda_{m}\in k$ for all 
$m\in \IN.$ Since $\rho_R(\p)$ is the restriction of $\rho_K(\p)$ to 
$V(R)$  (i) holds.

If  $\p(y) = 0$ then
$\lambda_{m} + y = (\text{\rm Tr}\, \rho_R(\p)^ m ) + y = \text{\rm 
Tr}\, \rho_{R/y}(\p(y))^ m = 0$ so $\lambda_{m} = 0.$ Thus $\text{\rm Tr}\, \rho_K(\p)^ 
m = 0$ for all $m$ and therefore $\rho_{K}(\p) = 0.$ Applying this to 
the adjoint representation yields that $\p = 0.$ This 
finishes the proof in the integral domain case.

\noi Step 3: $\fX$ {\it connected and reduced.}  For $x \in \fX$ we view $\p(x)$ 
as a $\kdiag$ element of $\fg(R/x)$ (see last Lemma). Fix $m \in \Bbb N$ and let $r = \text{\rm 
Tr}\, 
\rho_R(\p)^ m  \in R$ and $r_{x} = \text{\rm Tr}\, \rho_{R/x}(\p(x))^ m 
\in R/x.$ The integral 
domain case yields that $
r_{x} \in k \subset R/x.$ Clearly
$r_{y} = 
r_{x} + y$ whenever $x \subset y,$ and therefore $r_{x} \in k$ is 
constant on the irreducible components of $\fX$ (take $x$ to be a 
minimal ideal),  hence constant on $\fX$ (by [EGA] Ch.0 Cor.2.1.10 
since $\fX$ is connected and may be assumed noetherian). Call this common value $\lambda.$  
Since $\lambda = r_{x} = r + x$ for all $x$
and $R$ is 
reduced it follows that $r = \lambda$, hence that $r \in k.$   

Finally assume $\p(y) = 0.$ Let $x \subset y$ be a minimal prime, and 
view 
$\p(x)$ as a $\kdiag$ element of $\fg(R/x)$ and $y$ as an element of 
$\Spec (R/x)$. Since $\p(x)(y) = \p(y) = 0$,  the 
integral domain case yields that $\p(x) = 0$. Thus $\p$ vanishes 
in the irreducible component corresponding to $x$ and hence everywhere 
in $\fX$ as we saw above. Since $R$ is reduced this forces $\p = 0.$ \qed
\vs 
\noi {\bf Remark.} For a given $\lambda \in k$ let 
$\fg(\p)^{\lambda} = \{v \in \fg(R) : [\p,v] = \lambda v \}$. Then 
$\fg(\p)^{\lambda}$ is a projective $R$-submodule of $\fg(R)$. 
Conjugacy is related to the freeness of these submodules.

\vs

\demo {Proof of Theorem 1(i)} Let $\frak a$ be an abelian $\kdiag$ 
subalgebra of $\fg(R).$  Assume 
first that in addition of being connected  $\fX$ is reduced.  Then by 
Lemma 5(i) and Proposition 6(ii) elements of $\frak a$ are 
linearly independent if and only if they are so after evaluation at 
any element of $\fX.$ It follows that it will suffice to establish 
our result under the assumption  that $R$ is an integral 
domain.  Let then $\fa\subset \frak k$ and $K$ be as in Proposition 
4,
and fix a base $\Pi = \{\alpha _1\dots \alpha _\ell\}$  of the corresponding root
system $\Delta =\Delta (\fg(K),\frak k).$ If
$\omega ^\vee_1,\dots,\omega ^ \vee_\ell$ are the fundamental coweights corresponding to
the $\alpha _i\text{'s}$ then by reasoning as in Step 2 of the last 
Proposition we conclude that  $\fa \subset \oplus
k\omega ^\vee_i.$ This finishes the proof in the reduced case. 

 In general let $\J$ be the nilradical of $R$ and set $R' = R/\J.$ It 
 is clear that the image $\fa'$ of $\fa$ under the canonical map 
 $\fg(R) \rightarrow \fg(R')$ is abelian and $\kdiag$. It then 
 follows  from the reduced case that  for any given elements  
$\{\p_{1} 
\ldots \p_{\ell +1} \}$ of $\fa$, we can find  a nontrivial linear dependence relation 
$c_{1}\p_{1}' + 
\ldots + c_{\ell + 1}\p_{\ell +1}' = 0$ (where of course the $c_{i}'s$ 
depend on the elements, and $\p_{i}'$ 
denotes the image of $\p_{i}$ under the canonical map). Consider now the element  
$\p := c_{1}\p_{1}  + 
\ldots + c_{\ell + 1}\p_{\ell +1} \in \fa.$ Then the coordinates of 
$\p$ with respect to any basis of $\fg$ (viewed as an $R$-basis of 
$\fg(R)$) are in $\J$. From this it follows that $\ad_{\fg(R)}\p$ is  nilpotent. 
On the other hand since $\p \in 
\fa$ we also have that  $\ad_{\fg(R)}\p$ is diagonalizable (as a 
$k$-linear endomorphisms of $\fg(R))$. It 
follows that $\ad_{\fg(R)}\p = 0$ and hence that $\p = 0$ since 
$\fg(R)$ has trivial centre.\qed
\enddemo
\proclaim{Corollary}Let $\fh$ be a split Cartan subalgebra of $\fg$.  Assume 
$\fX = \Spec(R)$ is connected. Then $\fh$ is  the unique MAD of $\fg(R)$ contained in $\fh(R)$. 
\endproclaim
\demo{Proof} Clearly $\fh \subset \fg(R)$ is abelian and $\kdiag$, 
hence maximal because of its dimension. If $\fk \subset \fh(R)$ is an 
abelian $\kdiag$ subalgebra of $\fg(R)$,  then $\fh + \fk$ is also 
abelian and $\kdiag$. Again a dimension argument shows that $\fk 
\subset \fh$. \qed
\vs 
\noi {\bf Remark.} One can also give a direct proof of this Corollary. 
Note that the 
connectness assumption is also necessary. Indeed if $h \in \fh$ and 
$e \in R$ is an idempotent, then $h \otimes e$ is a $\kdiag$ element of $\fg(R)$ 
commuting with $\fh$.
\vs

The next result is crucial. Its effect is that the structure groups of 
the torsors related to conjugacy are {\it connected}.     
 
\proclaim {Proposition 7} Let $\fX = \Spec(R)$ be connected reduced and with a rational 
point. Let $\p \in \fg (R)$ be $k\text{-diagonalizable.}$ 
Fix  $x_0\in \fX$  such that $k(x_0) =k$ and set $\p_0 : = \p (x_0).$ If $x\in \fX$  then $\p (x)$
and $\p_0$ (viewed as two elements of $\fg (x)$) are conjugate under
$\fG(x).$
\endproclaim 

\demo {Proof}  By Lemma 5(i) and Proposition 
4 (with $R = k$) there
exists a split Cartan subalgebra $\frak h$ of $\fg$ with $\p _{0} \in \fh.$ Now
$\fh(x)$ is a split Cartan subalgebra of $\fg (x)$ and since any two such are
conjugate under $\fG(x)$ ([Bbk2] Ch~ 8 \S3.3 Cor. to Prop.~10) there is no loss
of generality in assuming that both $\p(x)$ and $\p _0$ belong to $\fh(x).$

Under this assumption, were $\p_0$ and $\p(x)$ not conjugate under 
$\fG(x),$ they
would be separated by  a polynomial function $f\in S\big(\fh(x)^*\big)$ 
which is invariant
under the Weyl group $W$ of $\big(\fg (x),\fh(x)\big)$ 
(ibid. Remarqu\' e \S5.2, and \S8.4 Lemma~6).  Now any
such $f$ is a linear combination of functions of the form $z \mapsto
\;\text{Tr}\; \rho_{k(x)}(z)^ m$ with  $\rho:\fg \to \fgl(V)$  a finite 
dimensional representation of $\fg$ (ibid \S8.2, Cor.~2).
But this is not possible. Indeed,  functoriality combined with Proposition 
6(i) yields

$\text{Tr}\; \rho_{k(x)}\big(\p (x)\big)^ m
= \big(\text{Tr}\;\rho_R(\p)^ m\big)(x)
= \big(\text{Tr}\; \rho_R(\p)^m\big)(x_0)
= \;\text{Tr}\;\rho_{k(x)}(\p_0)^ m. 
$
\qed
\enddemo 
\vs
\proclaim{Proposition 8} Let  $\p_{0}$  be a $\kdiag$ element of 
$\fg.$ Then.

\item {\rm (i)} $\fG(k)\cdot \p_0\subset \fg $ is
a Zariski closed set.

\item {\rm (ii)} Let $J\triangleleft S(\fg^*)$ be the defining ideal of
the affine variety $\fG(k)\cdot \p_0.$ Then the elements of  $J$
vanish on
$\fG(S)\cdot
\p_0$ for any $k\text{-algebra} \; S.$
\endproclaim 

\demo {Proof} (i) Let $\bar k$ be the algebraic closure of $k$. By [Brl]
Theorem~  9.2(ii) $\fG(\bar k)\cdot \p_0$ is a closed
subset of $\fg (\bar k).$ Since the Zariski topology of $\fg (\bar k)$
induces that of $\fg $ it will suffice to show that
$$
 \fG(\bar k)\cdot \p_0\cap \fg = \fG(k)\cdot \p_0
\tag "(8.1)"
$$
Let $\idealq = \idealP\cdot
\p_0\in \fg$ for some $\idealP\in \fG(\bar k).$  It is easy to see 
that   $\idealq$ is $\kdiag$.
We now use Proposition 4 (again with $R = k$) and conjugacy of split Cartans to see that 
to establish (8.1), we may assume  that  
$\idealq\in \fh$ where $\fh$ is some fixed split Cartan containing 
$\p_{0}.$  We then have
two elements $\p_0$ and $\idealq$ of $\fh$ which are conjugate under $\fG(\bar
k).$  A standard argument using the Bruhat decomposition of $\fG(\bar k)$ shows
that
$\idealq = w(\p_0)$ for some $w$ in the Weyl group $W$ of $\Delta(\fg , \fh)$.
Since $w$ is the restriction to $\fh$ of an element of $\fG(k)$ ([Bbk2]
 Ch.8 \S5 no.3 Remarqu\'e) (8.1) holds.  This finishes
the proof of (i).

(ii) The defining ideal of $\fG(k)\cdot \p_0$ is $J:=\{f\in S(\fg^*):f$
vanishes in $\fG(k)\cdot \p_0\}.$  First we assume that $S$ is an integral 
domain. In this case we establish (ii)
 by showing below  that $E:=\fG (k)\cdot \p_0$ is dense on $\fG
(F)\cdot \p_0$ where $ F$ is the algebraic closure of
the quotient field of $S.$

Let $T(E,\overline E) = \{\P\in \fG(F): \P\cdot E\subset \overline
E\}$ (hereafter $^{-}$ denotes Zariski closure).  This is a
closed subset of the affine variety $\fG(F)\cdot \p_0 \subset \fg
( F)$  ([Bor] Proposition 1.7 ii). Since $\fG(k)$ is dense in
$\fG(F)$ (ibid 18.3) we obtain
$$
\fG(k)\subset T(E,\overline E)\Longrightarrow \overline{\fG(k)} \subset \overline
{T(E,\overline E)} = T (E,\overline E) \Longrightarrow \fG(F)\subset
T(E,\overline E).
$$
In other words: $T(E,\overline E) = \fG( F).$  It follows 
then that in $\fg(F)$ we have
$$
\fG(F) \cdot \p_0 = \fG(F)\cdot \fG(k)\cdot \p_0 =
\fG(F)\cdot E \subset \overline E =\overline{\fG(k)\cdot \p_0}
\subset \overline{\fG(F)\cdot \p_0} = \fG(F)\cdot \p_0.
$$ 
Thus $ \fG(F)\cdot \p_0  =\overline
{\fG(k)\cdot \p_0}$ as desired. This establishes our result for 
integral domains.

Let $S$ now be arbitrary. Given that $k[\fG]$ is an integral 
domain the elements of $J$ vanish at $\text{id} \cdot \p_0$ where $\text{id} \in \fG(k[\fG]).$
Since any element of $\fG(S)$ is of the form $\phi(\text{id})$ for some arrow $\phi: k[\fG]
\rightarrow S$ the results holds for $S$ by functoriality. \qed
\enddemo 
\vs
 
\proclaim{Proposition 9} Let  $\p_{0}$ be a $\kdiag$ element of $\fg$, 
and let $\fL$ be its
isotropy group  {\rm(}i.e. 
$\fL(S) =\{\P\in \fG(S): \P\cdot\p _0 =\p_0\}$ for any $S$ 
in $\kalg$ {\rm ).}   Let $\fh$ be a 
fixed split Cartan subalgebra of $\fg$ 
containing $\p_{0}.$ Then there exists a base $\Pi =\{\alpha _1\dots
\alpha _\ell\}$ of $\Delta = \Delta(\fg,\fh)$ and a subset $I\subset \{1,\dots,\ell\}$
such that

\item {\rm (i)}  $\fL$  is the standard Levi subgroup 
corresponding to $I.$ In particular, $\fL$ is a  split
(connected)  reductive algebraic group. 

\item {\rm (ii)} The derived group $\fG_{I}$ of $\fL$ is of simply 
connected type.

\item {\rm (iii)} If condition \text{\rm (TLT)}  on $\X$  holds for  $\fG$, then it 
also holds for $\fG_{I}$

\endproclaim
 
\demo{Proof} Let $\Delta 
(\p_{0}) = \{ \alpha \in \Delta : \langle \alpha ,\p_{0}\rangle =0 \}.$ If $\Delta 
(\p_{0}) = \emptyset$ set $I =  \{1,\dots,\ell\}$. If not, then 
$\Delta (\p_{0}) $ is  a root system on 
the subspace of $\fh^{*}$ it spans,
and   there exists a subset $I = I(\p_{0})\subset \{1,\dots,\ell\}$
such that $\Pi_I: = \{\alpha _i:i\in I\}$ is a base of $\Delta (\p_{0}).$ ([Bbk1]
Ch.VI Proposition 24.)
 
\noi (i)   Over $\overline{k}$ this follows from Lemma 3.7 and Corollary 3.11 of [Stb]. 
For the general case one has to check that all  relevant arguments 
 hold over $k$ (e.g. the Bruhat decomposition of $\fG(k)$). 
 
 \noi (ii) By SGA3 Exp. XXII the derived group of $\fL$ is generated (as 
 a sheaf group on the flat site of $\fX$) 
 by the root subgroups corresponding to roots whose support lies in 
 $\Pi_{I}$. We denote this group by  $\fG_{I}$. That  
 $\fG_{I}$ is simply connected means that the geometric fibers of $\fG_{I}$ are 
 simply 
 connected algebraic groups, which holds by [Stb] 2.11 and 2,13 (See 
 also [SS] Ch. 2 Cor. 5.4). 
 
 \noi (iii)  We must show 
 that $\Hzar^{1}(\fX' , \fL') = (0)$ where $\fX' = \Xred$, and $\fL'$ is the standard Levi subgroup corresponding to a subset 
 $I'$ of $I$ (Strictly speaking our group is $\fL'_{\fX} := \fL' \times _{\Spec(k)} 
 \fX$ but we will omit the subindex $\fX$ for convenience. Similarly 
 for all the other  algebraic groups $\fG_{I}$, $\Gm$ etc. involved in 
 the proof). Consider the exact sequence (of sheafs of groups on the 
 flat site of $\fX$. See SGA3 Exp. XXII 6.2.3) 
 $$ 1 \rightarrow \fG_{I} \rightarrow \fL' \rightarrow \Gm^{r'} 
 \rightarrow 1.$$
  Since $\fL'$ is split, the above sequence 
splits and therefore it is exact on the Zariski site of $\fX$. 
Passing to C\v ech cohomology yields
$$0 \rightarrow \Hzar^{1}(\X', \, \fG_{I'}) \rightarrow 
\Hzar^{1}(\X', \, \fL' ) \rightarrow \Hzar^{1}(\X', \,  \Gm^{r'} ) 
\rightarrow 0.   $$
Now if condition (TLT) holds then $(0) = \Hzar^{1}(\X', \,  \fT) = 
\Hzar^{1}(\X', \,  \Gm)^{\ell} = \Pic(\X')^{\ell}$ and therefore 
$\Hzar^{1}(\X', \,  \Gm^{r'}) = \Pic(\X')^{r'} = (0)$. It follows that 
for establishing (iii), it suffices to show that $\Hzar^{1}(\X', \, 
\fG_{I'}) = (0)$. If now $\fL$ is the standard Levi subgroup of 
$\fG$ corresponding to $I'$, we can reason as above to conclude that 
the map $\Hzar^{1}(\X', \, \fG_{I'}) \rightarrow 
\Hzar^{1}(\X', \, \fL )$ has trivial kernel. Since under condition (TLT) 
we have
$\Hzar^{1}(\X', \, \fL ) = (0)$ the result follows \qed
 \enddemo
\vs
\noi {\bf Remark.} No argument is needed for Part (iii) 
of the last Proposition in the case of assumption (TRT) (since, 
unlike (TLT), this 
assumption
does not depend on $\fG$).
\comment
 
Passing to \'etale cohomology yields
 $$0 \rightarrow \Het^{1}(\X', \, \fG_{I'}) \rightarrow 
\Het^{1}(\X', \, \fL' ) \rightarrow \Het^{1}(\X', \,  \Gm^{r'} ) 
\rightarrow 0.   $$
Since the canonical map $\Hzar^{1}(\X' , \, -) \rightarrow \Het^{1}(\X' , \, -)$ 
is injective we obtain
$$0 \rightarrow \Hzar^{1}(\X', \, \fG_{I'}) \rightarrow 
\Hzar^{1}(\X', \, \fL' ) \rightarrow \Hzar^{1}(\X', \,  \Gm^{r'} ) 
\rightarrow 0.   $$
Now if condition (TLT) holds then $(0) = \Hzar^{1}(\X', \,  \fT) = 
\Hzar^{1}(\X', \,  \Gm)^{\ell} = \Pic(\X')^{\ell}$ and therefore 
$\Hzar^{1}(\X', \,  \Gm^{r'}) = \Pic(\X')^{r'} = (0)$. It follows that 
for establishing (iii), it suffices to show that $\Hzar^{1}(\X', \, 
\fG_{I'}) = (0)$. If now $\fL$ is the standard Levi subgroup of 
$\fG$ corresponding to $I'$, we can reason as above to conclude that 
the map $\Hzar^{1}(\X', \, \fG_{I'}) \rightarrow 
\Hzar^{1}(\X', \, \fL )$ has trivial kernel. Since under assumption (TLT) 
we have 
$\Hzar^{1}(\X', \, \fL ) = (0)$ the result follows \qed
 
\enddemo
\vs
\endcomment

\proclaim{Proposition 10} Let $\fX$, $\p$,  and $\p_{0}$ be as 
in Proposition 7. Let $J\triangleleft S(\fg^*)$ be the defining ideal of
the closed subset $\fG(k)\cdot \p_0 \in \fg,$ and $\fL$ the isotropy 
group of $\p_{0}.$ Then
 
\item {\rm (i)} There exists a canonical 
isomorphism $\fG/\fL \simeq \Spec (S(\fg^*)/J).$

\item {\rm (ii)} $\p$ vanishes on $J$ thereby inducing a scheme 
morphism $\psi_{\p}: \fX \rightarrow \fG/\fL.$
\endproclaim 
\demo{Proof}  
The abstract group $\fL(k)$ acts on the left on $k[\fG]$ via 
$\P \cdot f(\Q) = f(\P^{-1}\Q)$ for all $\P \in \fL(k), \, \Q \in 
\fG(k),$ and $f \in k[\fG]$ (where we are identifying $k[\fG]$ with 
the ring of polynomial functions of the Zariski closed set  
corresponding to $\fG(k)$). Since $\fL$ is reductive the quotient 
scheme $\fG/\fL$ exists and it is in fact  the affine scheme of the ring 
of invariants $ B := k[\fG]^{\fL(k)}$ ([MFK] Theorem1.1). There is a natural $k$-algebra 
homomorphism $S(\fg^*) \rightarrow B$ given by $\nu \mapsto f_{\nu}$ 
where $f_{\nu}(\Q) = \nu(\Q^{-1}\cdot \p_0).$ The kernel of this map is 
the defining ideal $J\triangleleft S(\fg^*)$  of
the closed set $\fG(k)\cdot \p_0.$ We thus have an injective 
$k$-algebra homomorphism $\phi : A \rightarrow B$ where $A =  
S(\fg^*)/J$. The surjectivity 
of $\phi$ follows from that of the induced homomorphism 
$\overline{\phi} := 1\otimes \phi : \bar k \otimes A \rightarrow \bar 
k \otimes B. $ Now to see that $\overline{\phi}$ is surjective (in 
fact an isomorphisms) it will suffice to show by [Brl] 9.1 that  $J$ generates 
in $S(\fg(\bar k)^{*})$ the defining ideal of $\fG(\bar k)\cdot 
\p_0.$ That $J$ has this property follows from Proposition 8(ii) 
applied to $S = \bar k.$

(ii)  We have for all $x \in \fX$
$$
S(\fg^*) \undersetbrace_{\p(x)}\to{\os \p\to\longrightarrow
R\os{\varepsilon _x}\to\longrightarrow} k(x).
$$
By Propositions 7 and 8(ii) $\p(x)$ vanishes on $J,$ thereby
inducing a homomorphism  $$
\overline \p:S(\fg^*)/J \simeq k[\fG]^\fL\to R.
$$
Indeed if
$f\in J$ then $f(\p):=\p(f)\in \us{x\in \fX}\to\cap x= (0) $ since $R$ is
reduced. Finally $\psi_{\p}$ is defined to be the scheme morphism corresponding to
$\overline\p.$ 

\enddemo

\vs\noi
\proclaim {Proposition 11} With the notation of Proposition 10 the
following are equivalent.

\item {\rm (i)} There exists $\P\in \fG(R)$ such that $\p_{0} =\P\cdot \p$
\item {\rm (ii)} There exists a scheme morphism $\widehat \psi _{\p} :\fX\to
\fG$ rendering the diagram
$$
\matrix\format \c &\c &\l\\
&&\fG\\
&\os \widehat \psi _\p\to\nearrow 
&\downarrow ^q\\ 
\fX&\underset{\psi _\p} \to\longrightarrow &\fG/\fL
\endmatrix
$$
commutative

\item {\rm (iii)} The pull back $\text{pr}_1: \fX \times _{\fG/\fL} \fG\to \fX$ admits a
global section.
\endproclaim 

\demo {Proof} It is  well known that (ii) and (iii) are equivalent.  To show
that (i) and (ii) are equivalent it is best to work in $k\text{-alg.}$
where by taking Proposition 10 into account the picture is as follows:

$$
\matrix\format\c &\c &\c\\
&&k[\fG]\\ 
&&\uparrow \\
&\os \idealP\to\swarrow &k[\fG]^\fL\\ 
&&\vert \wr \\
R &\os\overline\p\to\leftarrow &S(\fg^*)/J\\
&\us\p\to\nwarrow &\uparrow\\
&&S(\fg^*)
\endmatrix
$$

\noi Let $v_1,\dots,v_n$ be a basis of $\fg$ and $v^1,\dots ,v^n$ the
corresponding dual basis of $\fg^*.$  Then
$$
\p = \sum v_i\otimes \p(v^i).
$$
On the other hand (see the proof of Proposition 10 (i)) $ 
\idealP(v^i + J) =
v^i(\idealP^{-1}\cdot \p_0).$  In other words
$$
\idealP^{-1}\cdot \p_0 =\sum v_i\otimes \idealP(v^i + J).
$$
The commutativity of the diagram is thus equivalent to $\p$ and $\idealP^{-1}\cdot
\p_0$ being the same element.
\enddemo 
\qed

\vs\noi
{\bf Remark 12.} 
The picture that emerges after the
pull back by $\psi _\p$ is the following:
$$
\CD
\fX\times_{\fG/\fL} \fG @>{pr_2}>> \fG\\ 
@Vpr_1VV @VV{q}V\\ 
\fX @>>{\psi _\p}> \fG/\fL
\endCD
$$
Since the quotient morphism $q :\fG\to \fG/\fL$ is locally trivial (one can see this by means of the 
big cell, see [SGA]), the same is the
case for the   pullback $pr_1: \fX \times _{\fG/\fL} \fG\to \fX$
as a principal homogeneous space for $\fL$ over $\fX$ . Condition ({\bf TLT}) of Theorem 1(iii) ensures 
that $pr_1$ is trivial.
\vs
We now turn to the proofs of the last two parts of our main theorem.

\demo {Proof of Theorem 1(ii)(a) with $\fX$ reduced} Let ${\frak a}\subset \fg (R)$ be a regular 
MAD.  Fix a regular element
$\p \in \fa$.  Then for $\p_0$ as in Proposition 11 we have $\fL = 
\fT$ where $\fT$ is the split maximal torus of $\fG$ corresponding to 
a fixed split Cartan subalgebra $\fh$ of $\fg$ containing $\p_{0}$. Since
$\fT$ is a product of $\ell = \text{rank}(\fg)$ copies of the 
multiplicative group  $\Gm =\;\text{Spec}\; k[t^{\pm 1}],$ the
$\fT-$torsors over
$\fX$ are measured by
$$
\Het^ 1(\fX,\fT)\simeq \Het^ 1(\fX,\Gm)^\ell \simeq \;\text{Pic}\, (\fX)^ \ell.
$$
([Mln] Ch.4 \S4 and [DG] Ch.3 \S6.3]). The pull-back of Proposition 12(iii) is thus trivial and
we conclude that
$\p_{0} =\idealP\cdot \p$ for some $\idealP\in \fG(R).$ Then
$$
\idealP \cdot \fa \subset \idealP \cdot \frak z_{\fg (R)}\p =  
\frak z_{\fg(R)}\idealP \cdot\p = 
\frak z_{\fg(R)}\p_0 = \fh (R).
$$
Given that the only $k\text{-diagonalizable}$ elements of 
$\fh(R) $ are those of $\fh$ (see the Remark following the proof of Theorem 
1(i))  we have
$\idealP \cdot\fa\subset  \fh,$ and hence by maximality $\idealP \cdot 
\fa =
\fh$ as desired.
\enddemo

\demo {Proof of Theorem 1(ii)(b) with $\fX$ reduced}  By Proposition 11 
and Remark 12 if $\p \in \fa$ then $\P \cdot \p = \p_{0}$ for some $\P 
\in \fG(R)$. We may thus assume with no loss of
generality  that 
$\fa\cap \fg \neq (0).$ 
Fix a nonzero element $\p$ in this intersection as well as a split 
Cartan subalgebra
$\fh$ of $\fg$ with $\p\in \fh.$ We will reason by induction on the rank
$\ell$ of
$\fg.$ 

If $\ell =1$
then $\fg \simeq \frak{sl}_2$ so that $\p \ne 0$ amounts to $\p$ being 
regular and the result holds part (ii)(a). 
Assume now $\ell >1.$

Let $V_\p\subset \fh^*$ be the $k\text{-span}$ of those $\alpha \in \Delta
=\Delta (\fg,\fh)$ satisfying $\langle \alpha ,\p\rangle =0.$  If $\text{dim}_kV_\p =0$ then
$\p$ is regular.  We may thus assume
$0<\;\text{dim}_kV_\p <\;\ell .$  Let
$\Delta _\p =\Delta \cap V_\p.$  As mentioned in Proposition 9, $\Delta _\p$ is  a root system on $V_\p$
and   there exists a base $\Pi =\{\alpha _1\dots
\alpha _\ell\}$ of $\Delta $ and a subset $I\subset \{1,\dots,\ell\}$
such that $\Pi_I: = \{\alpha _i:i\in I\}$ is a base of $\Delta _\p.$ 

With such a $\Pi$ and $I$ fixed let 
 $\frak s$ be the subalgebra of $\fg$ generated by the $\fg ^{\pm \alpha
_i}\text{'s}$ with
$i\in I,$ and let $\fr
$ be the subalgebra of $\fg$ generated by $\fh$
and $\frak s.$   Then $\fr$ is reductive  with $\frak s$ as semisimple 
part, and centre $\frak c$ given by the $k\text{-span}$ of the coweights $\omega
^ \vee_j,$
$j\notin I.$ Note also that $\fr =\frak z_{\fg}(\p)$ and therefore 
that
$$
\fr(R) =\;\frak s(R)\oplus \frak c(R)
$$
 where $\frak c(R) $ is the centre  and
$\frak s(R)$ the derived algebra of $\fr(R).$ Since
$\p\in\frak a$ and
$\frak a$ is abelian we also have
$\fr(R) =\; \frak z_{\fg (R)}(\p) \supset \frak a.$

 Let
$$
\frak b =\{\ideals\in \frak s(R): \idealc + \ideals\in \frak a \ \ 
\text{some}\ \
\idealc\in
\frak c(R)\}.
$$
Then $\frak b =\pi (\frak a)$
where $\pi :\fr(R)\to \frak s(R)$ is the canonical homomorphism.
In particular $\frak a\subset \frak c(R) +\frak b$ and $\frak b$ is an abelian
$k\text{-diagonalizable}$ subalgebra of $\frak s(R).$  

Now $\frak b$ is
contained in a MAD of $\frak s(R)$. By induction  together with Proposition 
9(ii), we then deduce the existence of an element
$\idealP\in \fG_I(R)$  such that
$\idealP\cdot \frak b\subset \fh_\p :=\underset{\alpha \in\Delta _\p}\to\sum
[\fg^\alpha ,\fg^{-\alpha }].$  Since the elements of $\fG_I(R)$ fix $\frak
c(R)$ pointwise   we get
$$
\idealP\cdot \frak a\subset \idealP\cdot \big(\frak c(R)+\frak b\big) \subset
\frak c(R)+\fh\subset \fh (R).
$$
As before given that elements of $\idealP\cdot \frak a$ are $k\text{-diagonalizable}$
we in fact have $\idealP\cdot \frak a\subset \frak h$, and hence by 
maximality that
$\idealP\cdot
\frak a=\frak h$. 
\enddemo 
\demo{End of the proof of Theorem 1}
 Let $\J$ be 
the nilradical of $R$ and let $R' = R/\J.$ Then $\fXred = \Spec (R')$ is connected, has 
a rational point, and  by assumption satisfies property TLT. It follows 
then from 
Proposition 11 and Remark 
12 that if $\p'$ denotes the natural image of 
$\p$ in $\fg(R'),$ then 
$\P' \cdot \p' = \p'(x_{0} + \J) = \p_{0} \in \fg \subset 
\fg(R')$ for some $\P' \in \fG(R').$ We claim that there exists $\P \in \fG(R)$ lifting $\P'$ 
and such that $\P \cdot \p = \p_{0}$. This will finish the proof since we 
can then  reason as in the proofs of the reduced case above.

To establish the claim we may assume that $R$ is noetherian. In this 
case $\J$ is nilpotent and by considering $\J \supset \J^{2} \supset \J^{4} \ldots 
\supset \J^{2^{n}} = (0)$ it will suffice to establish the claim under the 
assumption  
$\J^{2} = (0).$
Since $\fG$ is smooth it then follows that $\P'$ does lift to an 
element $\P_{1}$ of $\fG(R)$. Thus
$$ \P_{1} \cdot \p = \p_{0} + \sum_{n=1}^{s}\alpha_{i}^{\vee} \otimes 
\epsilon_{i}+  
\underset{\alpha \in\Delta (\fg,\fh)}\to\sum v_{\alpha} \otimes \epsilon_{\alpha}$$
where $\{\alpha_{i}¥^{\vee}, v_{\alpha}\}$ is a Chevalley basis of 
$\fg$ and the $\epsilon_{i}$'s and $\epsilon_{\alpha}$'s belong to $\J$.

For  $\alpha \notin \Delta_{0} := \{\alpha \in  \Delta(\fg,\fh) : <\alpha , 
\p_{0}> = 0 \}$ let $\theta_{\alpha} = 
\exp(\ad( v_{\alpha} \otimes <\alpha,\p_{0}>^{-1}\epsilon_{\alpha})).$ This is an 
automorphisms of $\fg(R)$ that can be realized as the adjoint action 
of an element $\P_{\alpha}$ of $\fG(R)$ ([DG] II \S 6.3.7). If we now set $\P_{2} = 
\underset{\alpha \notin\Delta_{0}}\to\prod \P_{\alpha} $ (the 
product taken in any order) and $\P = \P_{2} \P_{1} \in \fG(R)$ we 
have $ \P \cdot \p = \p_{0} + \q $ where $\q = \underset{}\to\sum \alpha_{i}^{\vee} \otimes 
\epsilon_{i} + \underset{\alpha \in\Delta_{0}}\to\sum v_{\alpha} \otimes 
\epsilon_{\alpha}.$ Since $ \P \cdot \p$ is $\kdiag$ and commutes with $\p_{0}$ it follows that $\ad \, \q$ is 
$\kdiag.$ On the other hand $\ad \, \q$ is visibly nilpotent. Thus 
$\q = 0$ and $\P \cdot \p = \p_{0}$ as desired. \qed

\vs\noi
{\bf 13 An interesting example}
\vs\par\npb
We look at  $\fg=\frak {sl}_2$ and $\fX = \fG/\fT$ (the ``generic'' 
regular case for $\frak {sl}_2$).  
 Here the  group  
$\fG$ is $\text{\bf SL}_2$   and $R=k[\fG]^{\fT(k)}.$ For convenience we will denote $k[\fG]$ 
by $S.$ Consider the element
$\text{id} \in \fG(S).$  Let

$$
\p :=\;\text{id}^{-1}\cdot h =\pmatrix x_{22} &-x_{12}\\ \vvs -x_{21} &x_{11}\endpmatrix \;
\pmatrix 1 &0\\ \vvs 0 &-1\endpmatrix\;
 \pmatrix x_{11} &x_{12}\\ \vvs
x_{21} &x_{22}\endpmatrix=
$$

$$
\pmatrix x_{11}x_{22} + x_{12}x_{21} & 2x_{12}x_{22}\\ \vvs
-2x_{11}x_{21} &-(x_{11}x_{22} + x_{12}x_{21}) \endpmatrix\;.
$$

\noi Observe that $\p\in \fg(R).$  Furthermore $\p$ is regular (since $h$ is)
and $\kdiag$  (since $\text{ad}_{\fg(S)}\p$ is 
and  $\text{ad}_{\fg(R)}\p$ is simply its restriction to $\fg(R)$). 

 Say $\q \in \fg (R)$ is such
that $[\p,\q] = 2{\q}.$ Again by looking inside $\fg (S)$ we see
that
$$
\q =\;\text{id}^{-1}\cdot \pmatrix 0 &s\\ 0 &0\endpmatrix =
\pmatrix sx_{21}x_{22} &sx^2_{22}\\ \vvs
sx^ 2_{21} &sx_{21}x_{22}\endpmatrix\; .
$$
for some $ s\in S.$

Observe that, save for the $s,$ all entries of $\q$ belong to 
$S^ \alpha$. 
It follows  that $\q \in \fg(R)$ only if $s \in S^{-\alpha}.$ 
Since the $R$-module $S^{-\alpha}$ is (rank one projective but) not free we conclude that 
$\text{ad}_{\fg (R)}\p$  is
$k\text{-diagonalizable}$ as a $k\text{-linear}$  but {\it not} as an
$R\text{-linear}$ endomorphism of $\fg(R)$.

Next we look at the $k\text{-algebra}$ homomorphism $\overline{\p}$ attached to $\p$  
described in Proposition 10.   Let $\{E,H,F \}$ be the basis of $\fg^*$ dual 
to $\{e , h , f \}.$ Identify $S(\fg ^*)$ with the polynomial ring
$k[E, H, F].$  Then 
$ \p : S(\fg^*) \to R$ is given by
$$
\p =\cases
E &\mapsto 2x_{12}x_{22}\\ \vvs
H &\mapsto x_{11}x_{22} + x_{12}x_{21}\\ \vvs
F&\mapsto -2x_{11}x_{21} .
\endcases
$$
Our present situation is depicted by the diagram
$$
\matrix\format\c &\c &\c\\ \vvs
&&S\\ 
&&\uparrow\\
&&R\\
&&\vert \wr\\
R &\os \overline\p\to\leftarrow &S(\fg^*)/J\\
&\us\p\to\nwarrow &\uparrow\\
&&S(\fg^*)
\endmatrix
$$
and we can identify  $\overline{\p}$  with an endomorphism of $R.$ Our choice of $x_0 \in \fX$
is the maximal ideal of $R$ obtained by intersecting $R$ with the ideal of $S$ generated
by $x_{11} - 1$, $x_{12}$, $x_{21}$ ,and $ x_{22} - 1.$ Then $\p_0 = h$ and under our
isomorphism to $E + J \in
 S(\fg^*)/J$ corresponds $f_{E} \in R$ with 
$f_{E}(id) = E(id^{-1} \cdot \p_0) = 2x_{12}x_{22}.$ Similar considerations apply to $H$ and
$E$ thus showing that the endomorphism $\overline \p$ we are after is 
in fact the identity map.

According to Proposition 11 then, conjugacy is equivalent to the principal $\fT\text{-bundle}$
$q : \fG \rightarrow \fG/\fT$ being trivial. This however is not the case as the bundle in
question is a generator of $\text{Pic}(\fX) \simeq \Bbb Z.$ Note that $S/R$ is fppf and that our
bundle becomes trivial over $S$ as one can see directly (since by construction $\p$ is conjugate
to
$\p_0$ under $\fG(S)),$ or abstractly (since $\text{Pic}(\fG)$ is trivial).

\vs
\noi {\bf Acknowledgement.} I am thankful to S. Donkin and O. Mathieu 
for the interest they have shown on this work, as well as 
for their many helpful comments.

\vsk

\end
\Refs
\widestnumber\key{SGA3}

\ref\key {Bbk1}\by N. Bourbaki\book Groupes et alg\`ebres de Lie, Ch. 4,5, et 6.
\publ Hermann \yr 1968\endref

\ref\key {Bbk2}\by N. Bourbaki\book Groupes et alg\`ebres de Lie, Ch. 7 et 8
\publ Hermann \yr 1975\endref

\ref\key {Brl}\by A. Borel\book Linear algebraic groups. 2nd edition
\publ GTM 126 Springer Verlag \yr 1991\endref

\ref\key {CTO}\by J.-L. Colliot-Th\'el\`ene and M. Ojanguren\paper Espaces
principaux homog\`enes localement triviaux \jour Publ. IHES {\bf 72} \yr 1992
\pages 97-122\endref

\ref\key {CTS}\by J.-L. Colliot-Th\'el\`ene and J.-J. Sansuc\paper Fibr\'es quadratiques et
composantes connexes r\'eelles
\jour Math. Annalen {\bf 244}
\yr 1979\pages 105-134\endref

\ref\key {DG}\by M. Demazure et P. Gabriel\book Groupes
Alg\'ebriques\publ Tome I, North Holland\yr 1970\endref
 
\ref\key {Gll}\by P. Gille \book Open correspondence \yr 2001 \endref
 
\ref\key {EGA}\by  A. Grothendieck et J. A. Dieudonn\'e\book 
 El\'ements de G\'eometrie Alg\'ebrique
 \publ Grun. Math. 166. Springer Verlag \yr 1971 \endref

\ref\key {Jnz}\by C. Jantzen \book Representations of  algebraic groups
\publ Acad. Press \yr 1987\endref


\ref\key {Mln} \by J.S. Milne\book \'Etale cohomology\publ Princeton University
Press \yr 1980\endref


\ref\key {MFK}\by D. Munford, J. Fogarty and F. Kirwan\book Geometric 
invariant theory
\publ Springer-Verlag\yr 1994\endref

\ref\key {PK} \by D. Peterson and V. Kac\paper Infinite flag varieties 
and conjugacy theorems\jour
Proc. Natl. Acad. Sci., USA{\bf 80} \yr 1983\pages 1778-1782\endref


\ref\key {Pzl1}\by  A. Pianzola \paper Line bundles and conjugacy 
theorems for toroidal Lie algebras\jour 
 C.R. Acad. Sci. Canada {\bf 22} (3) \yr 
2000 \pages 125-128 \endref

\ref\key {Pzl2}\by  A. Pianzola \paper Affine Kac-Moody Lie algebras 
as torsors over the punctured line \jour 
Indagationes Mathematicae N.S. {\bf 13} (2)  \yr 2002 \pages 249-257 \endref

\ref\key {Pzl3}\by  A. Pianzola \paper Automorphisms of toroidal Lie algebras and their 
central quotients \jour Jour. of Algebra and Applications. {\bf 1}(1) 
\yr (2002) \pages 113-121 \endref

\ref\key {Rgn1}\by M. S. Raghunathan \paper Principal bundles on 
affine space and bundles on the projective line. \jour Math. Annalen {\bf 
285}
\yr 1989 \pages 309-332 \endref

\ref\key {Rgn2}\by M. S. Raghunathan  \paper Principal bundles on 
affine space  \jour Studies in Mathematics {\bf 8} \pages 187-206 Springer-Verlag \yr 1978 \endref

\ref\key {RR}\by M. S. Raghunathan and A. Ramanathan \paper Principal bundles on the
affine line \jour Proc. Indian Acad. Sci. {\bf 93}
\yr 1984 \pages137-145\endref

\ref\key {SGA3} \book Sch\'emas en Groupes (S\'eminaire dirig\'e par 
M. Demazure et A. Grothendieck) 
 \publ Lect. Notes Math. 151, 152, and 153. Springer Verlag \yr 1970 \endref

\ref\key {Slg}\by G. Seligman\book Rational methods in Lie algebras
\publ M. Dekker Lect. Notes in Pure \& Appl. Alg. Vol 17 \yr 1976 \endref

\ref\key {Stb}\by R. Steinberg  \paper Torsion in reductive groups  \jour 
Advances in Mathematics {\bf 15 }
 \pages 63-92  \yr 1975 \endref
 
 \ref\key {SS}\by  T. A. Springer and R. Steinberg  \paper Conjugacy 
 classes \jour In ``Seminar in algebraic groups and related finite 
 groups'' Lect. Notes in Math 131 \publ Springer Verlag
 \yr 1970\endref 
\endRefs
